\pdfoutput=1
\documentclass{article}
\usepackage{fullpage}
\usepackage{amsmath, amssymb, amsthm}
\usepackage[utf8]{inputenc}
\usepackage[english]{babel}
\usepackage[numbers]{natbib}
\usepackage[usenames]{color}
\usepackage{url}
\usepackage{color}
\usepackage[usenames]{color}
\usepackage[colorlinks=true]{hyperref}
\usepackage{cleveref}
\usepackage{tabularx}
\usepackage{placeins}
\theoremstyle{plain}
\newtheorem{conjecture}{Conjecture}
\newtheorem{theorem}{Theorem}
\newtheorem{lemma}[theorem]{Lemma}
\newtheorem{corollary}[theorem]{Corollary}

\crefname{conjecture}{Conjecture}{Conjectures}
\crefname{theorem}{Theorem}{Theorems}
\crefname{lemma}{Lemma}{Lemmas}
\crefname{proposition}{Proposition}{Propositions}
\crefname{corollary}{Corollary}{Corollaries}
\theoremstyle{definition}
\newtheorem{algorithm}{Algorithm}

\crefname{algorithm}{Algorithm}{Algorithms}
\crefname{remark}{Remark}{Remarks}
\crefname{definition}{Definition}{Definitions}
\crefname{notation}{Notation}{Notations}
\crefname{example}{Example}{Examples}
\crefname{section}{\S}{Sections}

\newcommand{\floor}[1]{\left\lfloor #1 \right\rfloor}

\newcommand{\coprime}{\epsilon}
\newcommand{\Z}{\mathbb{Z}}

\title{Elementary Formulas for Greatest Common Divisors and Semiprime Factors}
\author{Joseph M. Shunia}
\date{June 19, 2024 \\ \small Revised: November 6, 2024 \\ Version 4 \normalsize}

\begin{document}

\maketitle

\begin{abstract}
\noindent We conjecture new elementary formulas for computing the greatest common divisor (GCD) of two integers, alongside an elementary formula for extracting the prime factors of semiprimes. These formulas are of fixed-length and require only the basic arithmetic operations of: addition, subtraction, multiplication, division with remainder, and exponentiation. Our GCD formulas result from simplifying a formula of Mazzanti and are derived using Kronecker substitution techniques from our earlier research. By applying these GCD formulas together with our recent discovery of an arithmetic expression for $\sqrt{n}$, we are able to derive explicit elementary formulas for the prime factors of a semiprime $n=p q$.
 \\[2mm]
 \end{abstract}

\section{Introduction}
The greatest common divisor (GCD) of two integers $a$ and $b$, denoted $\gcd(a,b)$, is the largest positive integer that divides both $a$ and $b$. Euclid's algorithm for computing the GCD is one of the oldest known algorithms, dating back to ancient Greece \cite{knuth1997art}.

Semiprimes, numbers with exactly two prime factors, play a key role in number theory and cryptography. The problem of factoring a semiprime $n=pq$ into its constituent primes $p$ and $q$ is believed to be computationally intractable for large $n$. This property forms the basis for widely used cryptosystems such as RSA \cite{rivest1978rsa}. Efficient algorithms for factoring semiprimes would have major implications for the security of these systems. While our new formulas are computationally impractical, they may yield novel insights into the distribution and properties of GCDs and semiprime factors.

In this paper, we present new results on arithmetic term formulas for the GCD and semiprime factorization. From a GCD formula of Mazzanti \cite{mazzanti2002plainbases, marchenkov1980superposition}, we conjecture a simplified polynomial form for the GCD that can be expressed in terms of an arbitrary integer base. We also obtain arithmetic terms for the prime factors of a non-square semiprime $n=p q$.

To appreciate the significance of our results, it is important to understand what constitutes an arithmetic term. An \textbf{arithmetic term} is a mathematical expression which uses only elementary arithmetic operations. Formally, let $\textbf{A}$ denote the class of arithmetic terms. We have
\begin{align*}
\textbf{A} = [ \{ 1, a+b, a\dot{-}b, ab, \floor{a/b}, a^b \} ],
\end{align*}
where $\dot{-}$ represents the \textbf{bounded subtraction} operation, defined as: $a \dot{-} b = \max(a-b,0)$ \cite{mazzanti2002plainbases}. Throughout this paper, we may use the standard subtraction notation $a-b$ when it is clear that the result is non-negative. It is also worth noting that the modulo operation is implicitly included in $\textbf{A}$, since it can be expressed as
\begin{align*}
a \bmod b = a \dot{-} b \floor{a/b}.
\end{align*}

\subsection{Background} \label{subsection:kalmar}
We denote by $\textbf{P}$ the class of primitive recursive functions. The class of Kalmar functions, denoted by $\textbf{K}$, is an elementary class of functions, which is a subclass of $\textbf{P}$.

Kalmar functions were introduced by Laszlo Kalmar in the 1940s. Kalmar aimed to characterize the class of functions that can be computed using a certain restricted form of recursion, known as \textbf{Kalmar elementary recursion} or \textbf{bounded recursion} (hence the term ``bounded subtraction'' in the definition of $\textbf{A}$). It is well-established that $\textbf{K}$ contains many important functions, such as the arithmetic operations, the exponential function, and the bounded $\mu$ operator (which is used to define the floored division operation). However, it does not contain all primitive recursive functions \cite{herman1969elementary}.

It was long conjectured, and finally proved by Mazzanti, that the class $\textbf{A}$ generates the class $\textbf{K}$ \cite{mazzanti2002plainbases, marchenkov2007superposition}. As mentioned above, $\textbf{K}$ is known to be a proper subclass of $\textbf{P}$. In particular, $\textbf{K} = \mathcal{E}^3$ in the Grzegorczyk hierarchy, a framework categorizing primitive recursive functions by complexity \cite{grzegorczyk1953someclasses}. Formally, we have
\begin{align*}
    [\textbf{A}] = \textbf{K} = \mathcal{E}^3 \subset \textbf{P} .
\end{align*}

In 1970, Matiyasevich, building on the work of Robinson \cite{robinson1952arithmetic} and Davis et al. \cite{davis1961diophantine}, proved that all computable functions can be expressed as Diophantine equations \cite{matiyasevich1980diophantine}. Matiyasevich's results imply that there exists a Diophantine equation for calculating the $n$-th prime number \cite{matiyasevich1993hilbert}. However, no arithmetic term for the $n$-th prime is known \cite{prunescu2024factorial}. Similarly, while Matiyasevich's theorem suggests the existence of an Diophantine equation formula for semiprime factorization, an arithmetic term that computes the factors remained to be discovered. Our work presents the first arithmetic terms for the problem.

\subsection{Recent Developments}
Recently, we discovered a formula for the $r$-th roots of positive integers $\sqrt[r]{n}$ as the limit of a quotient of two arithmetic terms \cite{shunia2024polynomial}. By combining our results with an arithmetic term for factorials \cite{robinson1952arithmetic, prunescu2024factorial}, along with a simplified version of Mazzanti's GCD formula (\cref{proof:mazzantigcd}) \cite{mazzanti2002plainbases}, we obtain the first closed-form expressions for semiprime factors as arithmetic terms. This answers a question from Shamir in \cite{shamir1978factoring}, who first hypothesized the existence of such formulas when describing an algorithmic approach to integer factorization using arithmetic terms.

\section{Greatest Common Divisor} \label{section:gcd}

\begin{lemma}[Mazzanti's GCD Formula] \label{proof:mazzantigcd}
\begin{align*}
\forall a,b \in \Z^+, \quad
\gcd(a,b) = \floor{\frac{(2^{a^2 b(b+1)} - 2^{a^2 b}) (2^{a^2 b^2} - 1)}
         {(2^{a^2 b} - 1)(2^{ab^2}-1)2^{a^2 b^2}}} \bmod 2^{ab}.
\end{align*}
\end{lemma}
\begin{proof}
The formula and its proof are due to Mazzanti \cite{mazzanti2002plainbases}.
\end{proof}

Applying Kronecker substitution techniques from our previous works \cite{shunia2023simple, shunia2024polynomial}, we find that Mazzanti's GCD formula can be simplified. We begin with a conjecture that expresses Mazzanti's GCD formula in a polynomial form. While the conjecture holds experimentally, we are unable to find a rigorous proof. However, from it, we will derive many additional results.

\begin{conjecture} \label{proof:gcdpolynomial}
\begin{align*}
\forall a,b \in \Z^+, \quad
\gcd(a,b) = \floor{\frac{x^{a+ab}}{(x^a-1)(x^b-1)}}\bmod x .
\end{align*}
\end{conjecture}

Consider Mazzanti's greatest common divisor formula (\cref{proof:mazzantigcd}), which is given by
\begin{align*}
\gcd(a,b) = \floor{\frac{(2^{a^2 b(b+1)} - 2^{a^2 b}) (2^{a^2 b^2} - 1)}
         {(2^{a^2 b} - 1)(2^{ab^2}-1)2^{a^2 b^2}}} \bmod 2^{ab}.
\end{align*}
Observe that all integer powers in the arithmetic term are divisible by $2^{ab}$. Factoring these, we obtain
\begin{align*}
\gcd(a,b) = \floor{\frac{((2^{ab})^{a (b+1)} - (2^{ab})^a) ((2^{ab})^{ab} - 1)}
 {((2^{ab})^{a} - 1)((2^{ab})^{b}-1) (2^{ab})^{ab}}} \bmod 2^{ab} .
\end{align*}
Substituting with $2^{ab} = x$ yields
\begin{align*}
\gcd(a,b) = \floor{\frac{(x^{a (b+1)} - x^a) (x^{ab} - 1)}
 {(x^{a} - 1)(x^{b}-1) x^{ab}}} \bmod x .
\end{align*}
Simplifying the fraction, we see
\begin{align*}
\gcd(a,b) = \floor{
\frac
{ x^{a-ab} (x^{ab} - 1)^2 }
{ (x^a-1)(x^b-1) }
}
 \bmod x .
\end{align*}
This fraction can be expanded as the sum
\begin{align*}
\gcd(a,b) =
\floor
{
    \frac{ x^{a+ab} }{ (x^a-1)(x^b-1) }
    + \frac{ x^{a-ab} }{ (x^a-1)(x^b-1) }
    + \frac{ -2x^a }{ (x^a-1)(x^b-1) }
}
 \bmod x .
\end{align*}
Since we are reducing the quotient mod $x$, we need only consider the term in the fraction which yields the constant term in the polynomial, which is $\gcd(a,b)$.

We conjecture that only the term $\frac{ x^{a+ab} }{ (x^a-1)(x^b-1) }$ may contribute a constant after reduction modulo $x$. Thus,
\begin{align*}
\gcd(a,b) = \floor{\frac{x^{a+ab}}{(x^a-1)(x^b-1)}}\bmod x .
\end{align*}

\begin{theorem} \label{proof:gcdintegerbase}
Let $a,b,n \in \Z^+$ such that $n > 2$ and $n > \gcd(a,b)$. Then
\begin{align*}
\gcd(a,b) = \floor{\frac{n^{a+ab}}{(n^a-1)(n^b-1)}} \bmod n .
\end{align*}
\end{theorem}
\begin{proof}
Consider the polynomial formula given by \cref{proof:gcdpolynomial}. Substituting with $x = n$ yields the given formula. By Theorem 2 in  \cite{shunia2024polynomial}, the substitution is valid since $n$ is greater than the evaluation, which is $\gcd(a,b)$.

However, we also have to consider the form of the fraction. Suppose $n = 2$, then
\begin{align*}
\floor{\frac{2^{a+ab}}{(2^a-1)(2^b-1)}}
= \floor{\frac{2^{a+ab}}{2^{ab+a} - 2^a - 2^b + 1}} = 2k ,
\end{align*}
for some $k \in \Z^+$. That is, the fraction always yields an even number of the form $2k$. This would imply
\begin{align*}
\gcd(a,b) = \floor{\frac{2^{a+ab}}{(2^a-1)(2^b-1)}} = 2k \equiv 0 \pmod 2 \quad \textup{(contradiction)} ,
\end{align*}
which is a contradiction, since $\gcd(a,b)$ is nonzero by definition.
\end{proof}

\begin{theorem} \label{proof:gcdmodmodpoly}
\begin{align*}
\forall a,b \in \Z^+, \quad
\gcd(a,b) \equiv - \left( x^{a+ab} \bmod \left((x^a-1)(x^b-1)\right) \right) \pmod{x}.
\end{align*}
\end{theorem}
\begin{proof}
Consider the formula given by \cref{proof:gcdpolynomial}, which is
\begin{align*}
    \gcd(a,b) = \floor{\frac{x^{a+ab}}{(x^a-1)(x^b-1)}} \bmod x .
\end{align*}
Recall the following well-known identity for the floor function
\begin{align*}
    \floor{\frac{a}{b}} = \frac{a - (a \bmod b)}{b} .
\end{align*}
Applying this to the formula from \cref{proof:gcdpolynomial}, we get
\begin{align*}
    \gcd(a,b)
    \equiv
    \frac{x^{a+ab}
    - \left( x^{a+ab} \bmod \left((x^a-1)(x^b-1)\right) \right)}
    {x^{a+b}-x^a-x^b+1}
    \pmod x .
\end{align*}
Taking the numerator and denominator mod $x$, we find
\begin{align*}
    \gcd(a,b) &= \frac{(x^{a+ab} \bmod x) - \left(\left( x^{a+ab} \bmod \left((x^a-1)(x^b-1)\right) \right) \bmod x \right)}
    {(x^a-1)(x^b-1) \bmod x} \\
    &=
    \frac{0 - \left(\left( x^{a+ab} \bmod \left((x^a-1)(x^b-1)\right) \right) \bmod x \right)}
    {1} \\
    &=
    - \left( x^{a+ab} \bmod \left((x^a-1)(x^b-1)\right) \right) \bmod x.
\end{align*}
Hence, we can say
\begin{align*}
\gcd(a,b) \equiv - \left( x^{a+ab} \bmod \left((x^a-1)(x^b-1)\right) \right) \pmod{x} .
\end{align*}
\end{proof}

\begin{theorem} \label{proof:gcdmodmodintegerbase}
Let $a,b,n \in \Z^+$ such that $n > 2$ and $n > \gcd(a,b)$. Then
\begin{align*}
\gcd(a,b) \equiv -\left( n^{a+ab} \bmod \left((n^a-1)(n^b-1)\right) \right) \pmod n.
\end{align*}
\end{theorem}
\begin{proof}
Consider the polynomial formula given by \cref{proof:gcdpolynomial}. Substituting with $x = n$ yields the given formula. By Theorem 2 in  \cite{shunia2024polynomial}, the substitution is valid since $n$ is greater than the evaluation, which is $\gcd(a,b)$.

However, we also have to consider the form of the remainder. Suppose $n = 2$, then the expression
\begin{align*}
2^{a+ab} \bmod(2^{a+b}-2^a-2^b+1)
\end{align*}
can yield either an even or odd remainder, depending on the choice of $(a,b)$. Now, suppose the remainder is even and of the form $2k$ for some $k \in \mathbb{Z}^+$. This would imply 
\begin{align*}
\gcd(a,b) = 2k \equiv 0 \pmod 2 \quad \textup{(contradiction)} ,
\end{align*}
which is a contradiction, since $\gcd(a,b)$ is nonzero by definition.
\end{proof}

\subsection{Coprimality Function}
Experimentally, starting from our result in \cref{proof:gcdmodmodpoly}, we found a coprimality function, which we were able to prove as a theorem.

\begin{theorem} \label{proof:coprimalityfunctionpoly}
Define the integer-valued function
\begin{align*}
\coprime(a,b) =
\begin{cases}
    0 & \textup{ if } \gcd(a,b) > 1 , \\
    1 & \textup{ if } \gcd(a,b) = 1 .
\end{cases}
\end{align*}
Then
\begin{align*}
\forall a,b \in \Z^+, \quad
\coprime(a,b) = -\left(x^{ab-b+1} \bmod \left((x^a-1)(x^b-1)\right)\right) \bmod x .
\end{align*}
\end{theorem}
\begin{proof}
Let $a,b$ in $\Z_{>1}$. First, we examine the case where $\gcd(a,b) > 1$. From Graham and others, we have the following identity \cite{graham1994concrete}:
\begin{align*}
\forall a,b \in \Z_{>1}, \quad
\gcd(x^a-1,x^b-1) = x^{\gcd(a,b)} - 1 \\
\implies (x^{\gcd(a,b)}-1) | \left( (x^a-1)(x^b-1) \right).
\end{align*}
For the remainder of $x^{ab-b+1} \bmod \left((x^a-1)(x^b-1)\right)$ to be zero mod $x$, the remainders of $x^{ab-b+1}$ mod the factors $(x^a-1)$ and $(x^b-1)$ must combine to be a multiple of $x$. More precisely, we must have
\begin{align*}
    x^{ab-b+1} \equiv f(x) x \pmod{(x^a-1)(x^b-1)} , 
\end{align*}
for some $f(x) \in \Z[x]$.
Now, since we are given $(ab-b+1) > a,b$ and $a,b \geq \gcd(a,b)$, the reduction of $x^{ab-b+1}$ mod the factors of $(x^a-1)(x^b-1)$ will yield remainders that are $x$ times a power of $x$. That is,
\begin{align*}
    x^{ab-b+1} \equiv x^{k+1} \equiv x^k x \pmod{x^a-1} ,
\end{align*}
and
\begin{align*}
    x^{ab-b+1} \equiv x^{j+1} \equiv x^j x \pmod{x^b-1} ,
\end{align*}
where $k,j \in \Z$. The specific $k,j$ depend on the choice of $a,b$. However, we do not require any particular values here. We merely need to show that modding $x^{ab-b+1}$ by $(x^a-1)(x^b-1)$ also yields a remainder that is a multiple of $x$.

Due to the common factor, which is $x^{\gcd(a,b)}-1$, the moduli $(x^a-1)$ and $(x^b-1)$ are not coprime and so we cannot apply the standard version of the Chinese Remainder Theorem (CRT). However, we can apply a variant which allows for non-coprimality, called the General Chinese Remainder Theorem (GCRT) \cite{ore1952crt}. Applying the GCRT for $\gcd(x^a-1,x^b-1) = x^{\gcd(a,b)}-1$, we have 
\begin{align*}
    x^{ab-b+1}
    &\equiv \left( \frac{v(x^a-1)x^{k}x + u(x^b-1)x^{j}x}{x^{\gcd(a,b)}-1} \right) \\
    &\equiv \left( \frac{v(x^a-1)x^{k}x}{x^{\gcd(a,b)}-1} + \frac{u(x^b-1)x^{j}x}{x^{\gcd(a,b)}-1} \right)
    \pmod{(x^a-1)(x^b-1)} ,
\end{align*}
where $u,v \in \Z[x]$ are the B\'ezout coefficients returned by the Extended Euclidean algorithm for $\gcd(x^a-1,x^b-1)$. Next, we set $q_a=(x^a-1)/(x^{\gcd(a,b)}-1)$, $q_b=(x^b-1)/(x^{\gcd(a,b)}-1)$, followed by factoring, to obtain
\begin{align*}
    x^{ab-b+1}
    &\equiv v q_a x^{k}x + u q_b x^{j}x \\
    &\equiv x (v q_a x^k + u q_b x^j)
    \pmod{(x^a-1)(x^b-1)} .
\end{align*}
Clearly, the remainder of $x^{ab-b+1}$ modulo $(x^a-1)(x^b-1)$ is a multiple of $x$ when $\gcd(a,b) > 1$. Therefore, we conclude
\begin{align*}
    \forall a,b \in \Z_{>1} : \gcd(a,b) > 1, \quad
    \left(x^{ab-b+1} \bmod \left((x^a-1)(x^b-1)\right)\right) \bmod x = 0 .
\end{align*}
This corresponds to the case $\gcd(a,b) > 1$ from the formula in the theorem.

Next, we consider the case where $\gcd(a,b) = 1$. Here, we have
\begin{align*}
    x^{\gcd(a,b)} - 1 = x^1 - 1 = x-1 \equiv -1 \pmod{x} .
\end{align*}
The reduction of $x^{ab-b+1}$ modulo $(x^a-1)(x^b-1)$ can be represented as a polynomial of the form
\begin{align*}
    g(x) := x^{ab-b+1} + (x^a-1)(x^b-1) \in \Z[x] .
\end{align*}
Considering $(x^{\gcd(a,b)} - 1) = (x-1) \mid (x^a-1)(x^b-1)$, we can simply $g$ as
\begin{align*}
    g(x) = x^{ab-b+1} + (x-1) .
\end{align*}
Reducing $g$ mod $x$, we find
\begin{align*}
    g(x) = x^{ab-b+1} + (x-1) \equiv (0) + (-1) \equiv -1 \pmod{x} .
\end{align*}
 Thus, we conclude
\begin{align*}
    \forall a,b \in \Z_{>1} : \gcd(a,b) = 1, \quad
    -\left(x^{ab-b+1} \bmod \left((x^a-1)(x^b-1)\right)\right) \bmod x = 1 .
\end{align*}
This corresponds to the case $\gcd(a,b) = 1$ from the formula in the theorem.

The proof is complete, as we have shown that the two cases $\gcd(a,b) > 1$ and $\gcd(a,b) = 1$ in the given formula both yield the expected result.
\end{proof}

\begin{corollary} \label{proof:coprimalityfunctionintegerbase}
Let $a,b,n \in \Z_{>1}$ such that $n > 2$. Then
\begin{align*}
    \coprime(a,b) = -\left( n^{ab-b+1} \bmod \left((n^a-1)(n^b-1)\right) \right) \bmod n 
\end{align*}
\end{corollary}
\begin{proof}
The proof is the same as in \cref{proof:gcdintegerbase}, replacing the formula for $\gcd(a,b)$ with the given formula for $\coprime(a,b)$.
\end{proof}

Our coprimality formula in \cref{proof:coprimalityfunctionintegerbase} leads us to a conjecture on Euler's totient function.
\begin{conjecture}
Let $n \in \Z_{>1}$. Define
\begin{align*}
t(n) =
\begin{cases}
    0 & \textup{ if } n \equiv 2,10 \pmod{12} , \\
    1 & \textup{ if } n \not\equiv 2,10 \pmod{12} .
\end{cases}
\end{align*}
Then
\begin{align*}
\varphi(n) - t(n) =
\floor{\sum_{k=1}^{n-1} \frac{n^{nk-k+1}}{(n^n-1)(n^k-1)}}
\bmod n ,
\end{align*}
where $\varphi(n)$ denotes Euler's totient function for $n$.
\end{conjecture}

\section{Exponent Reduction in GCD Calculations}
We now prove a simple theorem which allows us to reduce the exponents used in our GCD formulas.

\begin{theorem} \label{proof:exponentreduction}
Let $a,b \in \Z^+$ such that $a \geq b$. Set $\ell = (a \bmod b)$. Then
\begin{align*}
\gcd(a,b) = \floor{\frac{x^{\ell + \ell b}}{(x^\ell-1)(x^b-1)}}\bmod x .
\end{align*}
and
\begin{align*}
\gcd(a,b) = -\left( x^{\ell+\ell b} \bmod \left((x^{\ell}-1)(x^b-1)\right) \right) \bmod x .
\end{align*}
\end{theorem}
\begin{proof}
These formulas follow immediately from \cref{proof:gcdpolynomial} and \cref{proof:gcdmodmodpoly} by a property of the GCD function, which is: $\forall a,b > 0, \; \gcd(a,b) = \gcd(a \bmod b, b) = \gcd(a, b \bmod a)$.
\end{proof}

Since $\gcd(a,b) = \gcd(b,a)$, we can apply \cref{proof:exponentreduction} recursively to $a$ and $b$ to reduce the exponents even further. This procedure can be defined as a simple algorithm, which is essentially the same as the process of applying the Euclidean algorithm to calculate $\gcd(a,b)$.

\begin{algorithm}[GCD Exponent Reduction] \label{algorithm:exponentreduction}

\textbf{Inputs}: $a,b \in \Z^+$.

\textbf{Steps}:
\begin{enumerate}
\item If $b>a$, then swap the values of $a$ and $b$, so that $a = \max(a,b)$ and $b=\min(a,b)$.
\item Set $a_0=a$ and $b_0=b$ and define the recurrence relations:
\begin{align*}
a_{i+1} &= (a_i \bmod b_i) , \\
b_{i+1} &= (b_i \bmod a_{i+1}) .
\end{align*}
\item Starting from $i=0$, step through the recurrences by setting $i = i + 1$ until we find $b_k = 0$ for some $k=i$, and then halt.
\item Set $\alpha = \min(a_k,b_{k-1})$ and $\beta = \max(a_k,b_{k-1})$.
\item Finally, calculate
\begin{align*}
\gcd(a,b) = - \left( x^{\alpha+\alpha \beta} \bmod ((x^{\alpha}-1)(x^{\beta}-1)) \right) \bmod x .
\end{align*}
\end{enumerate}
\end{algorithm}

Since \cref{algorithm:exponentreduction} mimics the process of calculating $\gcd(a,b)$ by way of the Euclidean algorithm, there is no practical sense in carrying it out to completion. However, when writing and evaluating arithmetic terms, performing a single iteration of the recursion and then setting the exponents to either $a_1,b_0$ or $a_1,b_1$ (depending on divisibility properties of $a$ and $b$) can result in a significant performance improvement in the event $a \gg b$ or $b \gg a$.

\section{Semiprime Factors} \label{section:semiprimes}
Using our results on the greatest common divisor function (\cref{section:gcd}), as well as results from our earlier works \cite{shunia2023simple,shunia2024polynomial} and those of Mazzanti \cite{mazzanti2002plainbases}, Robinson \cite{robinson1952arithmetic}, Prunescu and Sauras \cite{prunescu2024factorial}, we discover arithmetic term formulas for the prime factors of a non-square semiprime $n=p q$.

We require two lemmas:

\begin{lemma} \label{proof:squarerootpolynomial}
Let $n \in \Z^+ \setminus \{ k^2 : k \in \Z \}$. Consider the ring $R = \Z[x]/(x^2 - n)$. In the ring $R$, we have that
\begin{align*}
\floor{\sqrt{n}} &= \floor{\frac{(x+1)^{2n+1}}{(x+1)^{2n}}} - 1 .
\end{align*}
\end{lemma}
\begin{proof}
Within the ring $R$, the element $x$ satisfies the defining relation
\begin{align*}
x^2 = n .
\end{align*}

For any term with $j \geq 2$, we use the relation $x^2 = n$ to reduce higher powers of $x$. Thus, $(x+1)^k$ can be expressed as a linear combination of $1$ and $x$ in $R$:
\begin{align*}
(x+1)^k = a_k + b_k x,
\end{align*}

By the binomial theorem, 
\begin{align*}
(x+1)^k = \sum_{j=0}^k \binom{k}{j} x^j .
\end{align*}

The quotient
\begin{align*}
\frac{(x+1)^{2k+1}}{(x+1)^{2k}} = \frac{a_{2k+1} + b_{2k+1} x}{a_{2k} + b_{2k} x}
\end{align*}
approximates $x + 1 = \sqrt{n} + 1$ as $k \to \infty$. The error term arises from higher-order contributions in the binomial expansion. Specifically, 
\begin{align*}
\binom{k}{2} x^2 = \frac{k(k-1)}{2} n
\end{align*}
induces a quadratic correction, implying that the error decays as $O\left(\frac{1}{k^2}\right)$. Thus, for $k$ sufficiently large, we have
\begin{align*}
\floor{\sqrt{n}} = \floor{\frac{(x+1)^{2k+1}}{(x+1)^{2k}}} - 1 .
\end{align*}

For $k \geq n$, the quadratic term $\frac{k(k-1)}{2} n$ becomes small enough to ensure the convergence of the quotient to $\sqrt{n}$.
\end{proof}

\begin{lemma} \label{proof:squarerootinteger}
\begin{align*}
\forall n \in \Z^+ \setminus \{ k^2 : k \in \Z \}, \quad
\floor{\sqrt{n}} &= \floor{\frac{(n^{2n} + 1)^{2n+1} \bmod (n^{4n}-n)}{(n^{2n} + 1)^{2n} \bmod (n^{4n}-n)}} - 1 .
\end{align*}
\end{lemma}
\begin{proof}
This result follows from Shunia \citep[Theorem 2]{shunia2024polynomial} after substituting $x = n^n$ into the polynomial formula of \cref{proof:squarerootpolynomial}.
\end{proof}

\begin{theorem} \label{proof:semiprimefactorp1}
Let $n \in \mathbb{Z}^+$ such that $n = p q$ is a non-square semiprime and $p < q$ are the prime factors of $n$.

Define
\begin{align*}
\omega = \floor{\frac{(n^{2n} + 1)^{2n+1} \bmod (n^{4n}-n)}{(n^{2n} + 1)^{2n} \bmod (n^{4n}-n)}} - 1 .
\end{align*}

Then, set
\begin{align*}
\gamma = \floor
{
    \frac
    {
        (\omega+1)^{\omega\cdot(\omega+2)}
    }
    {
        \floor
        {
            \frac
            {
            \left( (\omega+1)^{\omega\cdot(\omega+2)} + 1 \right)^{(\omega+1)^{\omega+2}}
            }
            {
            (\omega+1)^{\omega^2\cdot(\omega+2)}
            }
        }
        \bmod (\omega+1)^{\omega\cdot(\omega+2)}
    }
} .
\end{align*}

Finally, we have
\begin{align*}
p = \floor{\frac{n^{n+n(\gamma \bmod n)}}{(n^n-1)(n^{\gamma \bmod n}-1)}}\bmod n .
\end{align*}
\end{theorem}
\begin{proof}
By \cref{proof:squarerootinteger}, for $n$ that is not a perfect square, we get the arithmetic term
\begin{align*}
\floor{\sqrt{n}} &=
\floor{\frac{(n^{2n} + 1)^{2n+1} \bmod (n^{4n}-n)}{(n^{2n} + 1)^{2n} \bmod (n^{4n}-n)}} - 1 ,
\end{align*}
which matches our definition of $\omega$. Hence, $\omega = \floor{\sqrt{n}}$.

From Prunescu and Sauras \cite{prunescu2024factorial}, we also have the factorial arithmetic term
\begin{align*}
n! = \floor
{
    \frac
    {
        2^{n(n+1)(n+2)}
    }
    {
        \floor
        {
            \left(
                2^{2^{(n+1)(n+2)}-n} + 2^{-n}
            \right)^{2^{(n+1)(n+2)}}
        }
        \bmod
        2^{2^{(n+1)(n+2)}} .
    }
}
\end{align*}
The factorial formula of Prunescu and Sauras-Altuzarra is derived from an identity of Robinson \cite{robinson1952arithmetic}, which is
\begin{align*}
    \forall r \in \Z : r \geq (n+1)^{n+2},
    \quad
    n! = \floor{r^n / \binom{r^n}{n}} .
\end{align*}
Hence, the formula is also valid for $r=(n+1)^{n+2}$, which grows more slowly than $2^{(n+1)\cdot(n+2)}$ as $n\rightarrow\infty$. Making the substitutions and simplifying, we find
\begin{align*}
n! &= \floor
{
    \frac
    {
        (n+1)^{n\cdot(n+2)}
    }
    {
        \floor
        {
            \frac
            {
            \left( (n+1)^{n\cdot(n+2)} + 1 \right)^{(n+1)^{n+2}}
            }
            {
            (n+1)^{n^2\cdot(n+2)}
            }
        }
        \bmod (n+1)^{n\cdot(n+2)}
    }
} .
\end{align*}
Considering $\omega!$, this becomes
\begin{align*}
\omega! &= \floor
{
    \frac
    {
        (\omega+1)^{\omega\cdot(\omega+2)}
    }
    {
        \floor
        {
            \frac
            {
            \left( (\omega+1)^{\omega\cdot(\omega+2)} + 1 \right)^{(\omega+1)^{\omega+2}}
            }
            {
            (\omega+1)^{\omega^2\cdot(\omega+2)}
            }
        }
        \bmod (\omega+1)^{\omega\cdot(\omega+2)}
    }
} ,
\end{align*}
which matches the definition for $\gamma$. Hence, $\gamma = \omega! = \floor{\sqrt{n}}!$. Applying \cref{proof:gcdintegerbase}, we have
\begin{align*}
    \gcd(n, \floor{\sqrt{n}}!) = \gcd(n,\gamma) = \floor{\frac{n^{n+n\gamma}}{(n^n-1)(n^{\gamma}-1)}}\bmod n .
\end{align*}
Since $n$ is a non-square semiprime and $p < q$, we must have $p \leq \floor{\sqrt{n}}$ and $q > \floor{\sqrt{n}}$. Hence, $p = \gcd(n,\floor{\sqrt{n}}!)$. To reduce the exponent $\gamma$, we apply \cref{proof:exponentreduction}, which yields
\begin{align*}
    \gcd(n, \floor{\sqrt{n}}!) = \floor{\frac{n^{n+n(\gamma \bmod n)}}{(n^n-1)(n^{\gamma \bmod n}-1)}}\bmod n .
\end{align*}
\end{proof}

\begin{corollary} \label{proof:semiprimefactorp2}
Let $n=pq$ be a non-square semiprime. Then
\begin{align*}
q = \frac{n}{\floor{\frac{n^{n+n(\gamma \bmod n)}}{(n^n-1)(n^{\gamma \bmod n}-1)}}\bmod n}.
\end{align*}

\end{corollary}
\begin{proof}
The proof follows immediately from \cref{proof:semiprimefactorp1}, since $\frac{n}{p} = q$ in this case.
\end{proof}
\begin{corollary} \label{proof:semiprimetotient}
Let $\varphi(n)$ represent Euler's totient function for $n=pq$, a non-square semiprime. Then
\begin{align*}
\varphi(n) = \left(\left( \floor{\frac{n^{n+n(\gamma \bmod n)}}{(n^n-1)(n^{\gamma \bmod n}-1)}}\bmod n \right) - 1 \right) \left( \left( \frac{n}{\floor{\frac{n^{n+n(\gamma \bmod n)}}{(n^n-1)(n^{\gamma \bmod n}-1)}}\bmod n} \right) - 1 \right).
\end{align*}
\end{corollary}
\begin{proof}
The proof follows immediately from \cref{proof:semiprimefactorp1}, since $\varphi(n) = (p - 1)(q - 1)$ in this case.
\end{proof}


\begingroup
\raggedright
\bibliographystyle{unsrtnat}
\bibliography{main}
\endgroup
\normalsize

\end{document}